\magnification=\magstep1
\input amstex
\documentstyle{amsppt}
\redefine\R{\Bbb R}

\redefine\P{\text{\rm Power}}
\redefine\B{\text{\rm Borel}}

\define\dom{\text{\rm dom}}
\define\rng{\text{\rm rng}}
\define\ba{\Bbb B}
\pageheight{7in}
\topmatter
\title
Forcing with ideals generated by closed sets
\endtitle
\author
Jind\v rich Zapletal
\endauthor
\affil
University of Florida
\endaffil
\abstract
Consider the posets $P_I=\B(\R)\setminus I$ where $I$ is a $\sigma$-ideal $\sigma$-generated by 
a projective collection of closed sets. Then the $P_I$ extension is given by a single real $r$ of an almost minimal degree: every real
$s\in V[r]$ is Cohen-generic over $V$ or $V[s]=V[r].$ 
\endabstract
\thanks
The author is partially supported by grants GA \v CR 201-00-1466 and NSF DMS-0071437.
\endthanks
\address
Department of Mathematics, University of Florida, Gainesville FL 32611
\endaddress
\email
zapletal\@math.ufl.edu
\endemail
\endtopmatter

\document

\subhead {0. Introduction}\endsubhead

Under suitable large cardinal hypotheses, every proper definable forcing for adding a single real is forcing-equivalent to the poset $\B(\R)\setminus I$
ordered by inclusion, 
for a suitable $\sigma$-ideal $I$ \cite {Z2}. In this paper I will analyze the case of $\sigma$-ideals $I$ $\sigma$-generated by 
a projective collection of closed sets.
For such an ideal the forcing $\B(\R)\setminus I$ is always proper. The representatives include some familiar posets (Sacks real=$\B(\R)$
minus the countable sets, Miller real=$\B(\omega^\omega)$ minus the modulo finite bounded sets, Cohen real=$\B(\R)$ minus the meager sets) as
well as posets as yet not used nor understood. Consider the forcing $\B(\R)\setminus$ the ideal of closed measure zero sets. Or, consider the ideal $J$ on $P(\omega)$ generated by the sets $A_x=\{y\subset\omega:y\subset x\}$ as $x$ varies over
all subsets of $\omega$ of asymptotic density $0$. The forcing Borel(asymptotic density $0$ ideal)$\setminus J$ naturally increases the cofinality of the
asymptotic density $0$ ideal.

The main result of this paper is

\proclaim {0.1. Theorem}
(ZFC+large cardinals) Let $I$ be a $\sigma$-ideal $\sigma$-generated 
by a projective collection of closed sets. The poset $P_I=\B(\R)\setminus I$ is proper and adds a single real $r_{gen}$ of
an almost minimal degree: If $V\subseteq V[s]\subseteq V[r_{gen}]$ is an intermediate model for some real $s$, then $V[s]$ is a Cohen
extension of $V$ or else $V[s]=V[r_{gen}].$
\endproclaim

\noindent However, the point of the paper is not exactly to prove this theorem. Rather, the point is to expose certain technologies that connect the 
descriptive set theory
with the practice of definable proper forcing. Another point is to show that there are certain posets about which one can prove quite a bit by virtue
of the syntax of their definition, but (to date) one can prove absolutely nothing from the combinatorics of the related objects.

The notation in the paper follows the set theoretic standard of \cite {J}. AD denotes the Axiom of Determinacy. If $I$ is a $\sigma$-ideal on the reals
the symbol $P_I$ denotes the poset $\B(\R)\setminus I$ ordered by inclusion. The real line $\R$ is construed to be the set of all total
functions from $\omega$ to $\omega.$ The large cardinal hypothesis needed for the proof of Theorem 0.1 can be specified to be ``infinitely many Woodin cardinals'' 
or less, depending on the descriptive complexity of the ideal $I.$

\subhead {1. The results} 
\endsubhead

\proclaim {1.1. Lemma}
If $I$ is a $\sigma$-ideal on the real line then $P_I$ forces ``for some unique real $r$ the generic filter is just
the set $\{B^V: B$ is a Borel ground model coded set of reals with $r\in B\}$''.
\endproclaim

\demo {Proof}
Let $\dot r$ be the $P_I$-name for a real defined by $\dot r(\check n)\in\check m$ if the set $\{s\in\R:s(n)=m\}$ belongs to the generic filter.
Note that this indeed defines a name for a total function from $\omega$ to $\omega$ since the collections $\{\{s\in\R:s(n)=m\}:m\in\omega\}$
are maximal antichains in the poset $P_I$ for each integer $n\in\omega.$ Also note that if any real is to be in the intersection of all
Borel sets in the generic filter, it must be $\dot r.$ In order to complete the proof, it is enough to argue by induction on the complexity
of the Borel set $B\notin I$ that $B\Vdash\dot r\in\dot B$ where $\dot B$ is the Borel set in the extension with the same Borel
definition as $B.$

Now this is clearly true if $B$ is a closed or a basic open set. Suppose $B=\bigcap_n A_n$ and for all $n\in\omega,$ $A_n\Vdash\dot r\in\dot A_n$
has been proved. Then for all $n\in\omega$ $B\subset A_n,$ so $B\Vdash\forall n\in\omega\ \dot r\in\dot A_n$ and $B\Vdash\dot r\in\bigcap_n\dot A_n=
\dot B$ as desired. And suppose that $B=\bigcup_nA_n$ and for all $n\in\omega,$ $A_n\Vdash\dot r\in\dot A_n$
has been proved. The collection $\{A_n:A_n\notin I\}$ is a predense set in $P_I$ below $B$--this simple observation uses the $\sigma$-completeness
of the ideal $I.$ So $B\Vdash\exists n\in\omega\ \dot r\in\dot A_n,$ and $B\Vdash\dot r\in\bigcup_n\dot A_n=\dot B$ as desired. Since
all Borel sets are obtained from closed sets and basic open sets by iterated applications of a countable intersection and union, the proof
is complete. \qed
\enddemo

The unique real from the statement of the previous lemma will be called the $P_I$-generic real and denoted by $\dot r_{gen}.$ This real is forced to
fall out of all ground model Borel sets in the ideal $I.$ Another standard piece of terminology: if $M$ is an elementary submodel of some large
structure and $r$ is a real such that the set $\{B\in P_I\cap M:r\in B\}$ is an $M$-generic filter on $P_I$, then I will call the real
$r$ $M$-generic.

\proclaim {1.2. Lemma}
If $I$ is a $\sigma$-ideal $\sigma$-generated by a closed sets then the poset $P_I$ is $<\omega_1$-proper.
\endproclaim

\demo {Proof}
Let me first show that the poset $P_I$ is proper. Suppose $A\in P_I$ is an arbitrary condition and $M$ is a countable elementary submodel
of a large enough structure containing all the relevant information. I must produce a master condition $A^*\subset A$ for the model $M.$
Consider the set $A^*$ of all $M$-generic reals in $A$. This set is Borel, and
if it is $I$-positive, by the previous lemma it forces $r_{gen}\in\dot A^*$, which is to say ``$\dot G\cap \check M\subset \check P_I\cap\check M$
is $\check M$-generic'', which is to say that $A^*$ is a master condition for the model $M$. Note also
that if $B\subset A$ is any other master condition for the model $M,$ necessarily $B\setminus A^*\in I.$ So $A^*$ is really the only candidate
for the required master condition. The only thing left to verify is $A^*\notin I.$

Suppose that $\{C_n:n\in\omega\}$ is a collection of closed sets in the ideal $I.$ I must produce a real $r\in A^*\setminus\bigcup_n C_n.$ Let
$D_n:n\in\omega$ be a list of all open dense subsets of the poset $P_I$ in the model $M,$ and by induction on $n\in\omega$ build
sets $A=A_0\supset A_1\supset A_2\supset\dots$ in the model $M$ so that for every $n\in\omega,$ $A_{n+1}\in D_n$ and $A_{n+1}\cap C_n=0.$
To perform the inductive step, first choose a set $A_{n+0.5}\subset A_n$ in $M\cap D_n$ and then note that since the set $A_{n+0.5}\setminus C_n$
is $I$-positive and the set $C_n$ is closed, there must be a basic open neigborhood $O_n$ such that $O_n\cap C_n=0$ and $A_{n+0.5}\cap O_n$
is still $I$-positive. But then the set $A_n=A_{n+0.5}\cap O_n$ is in the model $M$ and satisfies the inductive assumptions. Once the induction
is complete, look at the $M$-generic filter on $P_I\cap M$ generated by the sequence of $A_n$'s. By the previous lemma applied in the model $M,$
the intersection of all the sets in this filter is a singleton containing a real $r.$ By the construction, $r\in A^*\setminus\bigcup_n C_n$ as desired.

The attentive reader will have noticed that the previous argument gives even strong properness of the poset $P_I$, see \cite {S}. A slight variation of
the argument will give $<\omega_1$-properness.

By induction on $\alpha\in\omega_1$ prove that the poset $P_I$ is $\alpha$-proper. The successor step is trivial on the account of the previously proved
properness of $P_I$. So suppose that $\alpha$ is a limit ordinal, a limit of an increasing sequence $\alpha_0\in\alpha_1\in\alpha_2\in\dots$
and for all $n\in\omega$ the $\alpha_n$-properness of the poset $P_I$ has been proved. Let $A\in P_I$ be an arbitrary condition and let 
$\langle M_\beta:\beta\in\alpha\rangle$ be a continuous $\in$-tower of countable elementary submodels of a large enough structure
such that $A\in M_0.$ As before, it is enough to show that the set $A^*=\{r\in A:$ for all $\beta\in\alpha$ the real $r$ is $M_\beta$-generic
$\}$ is $I$-positive, since it will be the required master condition for the tower.

Suppose that $\{C_n:n\in\omega\}$ is a collection of closed sets in the ideal $I.$ I must produce a real $r\in A^*\setminus\bigcup_n C_n.$ Let $M=
\bigcup_{\beta\in\alpha}M_\beta$ and let
$D_n:n\in\omega$ be a list of all open dense subsets of the poset $P_I$ in the model $M$ such that $D_n\in M_{\alpha_n+1}$,
and by induction on $n\in\omega$ build
sets $A=A_0\supset A_1\supset A_2\supset\dots$ so that for every $n\in\omega,$ $A_{n+1}\in D_n\cap M_{\alpha_n+1},$ 
$A_{n+1}\subset\{r\in A:$ the set $\{B\in P_I\cap M_\beta:r\in B\}$ is an $M_\beta$-generic
filter, for all $\beta\in\alpha_n\}$, and $A_{n+1}\cap C_n=0.$
To perform the inductive step, first look at the set $A^*_n=\{r\in A_n:r$ is $M_\beta$-generic
for all $\alpha_{n-1}\in\beta\in\alpha_n\}$. This set is Borel, it is $I$-positive by the $\alpha_n$-properness of the poset $P_I$ (it is the only candidate
for the master condition for the tower $\langle M_\beta:\alpha_{n-1}\in\beta\in\alpha_n\rangle$) and it is in the model $M_{\alpha_n+1}$.
As in the second paragraph of this proof it is now possible to choose a set $A_{n+1}\subset A^*_n$ in $D_n\cap M_{\alpha_n+1}$ with $A_{n+1}\cap
C_n=0.$ Such a set satisfies the inductive assumptions. Once the induction
is complete, look at the $M$-generic filter on $P_I\cap M$ generated by the sequence of $A_n$'s. By the previous lemma applied in the model $M,$
the intersection of all the sets in this filter is a singleton containing a real $r.$ By the construction, $r\in A^*\setminus\bigcup_n C_n$ as desired.
\qed
\enddemo
 
It is well known that if $I$ is a $\sigma$-ideal such that the forcing $P_I$ is proper, and $B\Vdash\dot s$ is a real, then by using a stronger condition
$C\subset B$ if necessary the name $\dot s$ can be reduced to a Borel function $f:C\to\R$ such that $C\Vdash\dot s=\dot f(\dot r_{gen}).$ To see
how this can be done, choose a countable elementary submodel $M$ of a large enough structure, let $C$ be the set of all $M$-generic reals
in the set $B$ and define $f:C\to\R$ by $f(r)=\{\langle n,m\rangle:$ for some set $a\in P_I\cap M,$ $r\in A$ and $A\Vdash\dot s(\check n)=\check m\}$.
An absoluteness argument just like in the proof of the previous lemma shows that this function will work.

The following theorem is assembled from results of Martin and Solecki and appears in \cite {S}.

\proclaim {1.3. Lemma} 
(ZFC+large cardinals) If $I$ is a $\sigma$-ideal generated by closed sets, and if $A\subset\R$ is an $I$-positive 
projective set of reals then $A$ has a Borel $I$-positive subset.
\endproclaim

And the key tool for establishing Theorem 0.1 is

\proclaim {1.4. Lemma}
(projective uniformization) Suppose  $I$ is a $\sigma$-ideal such that

\roster
\item $I$ is generated by a projective collection of projective sets
\item every projective $I$-positive set has an $I$-positive Borel subset
\item the forcing $P_I$ is $<\omega_1$-proper in all forcing extensions.
\endroster

If $G\subset P_I$ is a generic filter and $V\subseteq V[H]\subseteq V[G]$ is an intermediate model, then $V[H]$ is a c.c.c.
extension of $V$ or else $V[G]=V[H].$
\endproclaim

Here (1) means that there is an integer $n$ such that the ideal is generated by boldface $\Sigma^1_n$ sets and the set of all codes for boldface 
$\Sigma^1_n$ sets in the ideal is itself projective.

\demo {Proof}
Let $I$ be a $\sigma$-ideal satisfying the assumptions of the lemma. On the account of (3) I can assume that the continuum hypothesis holds,
because it can be forced by a $\sigma$-closed notion of forcing, not changing the poset $P_I.$ Suppose that $\ba$ is a nowhere c.c.c. complete subalgebra of the completion
of the poset $P_I$. I will prove that the generic real $\dot r$ for the poset $P_I$ can be recovered from the generic filter 
$\dot H\subset \ba$. 

First, a piece of notation: Suppose $M$ is a countable elementary submodel of a large enough structure. I will say that a real $r$ is $M,\ba$-weakly
generic if for every $b\in \ba\cap M$ either there is a set $A\in P_I\cap M$ with $A\leq b$ and $r\in A$, or there is a set
$A\in P_I\cap M$ with $A\leq \lnot b$ and $r\in A.$ In such a situation, I will write $\dot H\cap M/r$ to denote the set $\{b\in\ba\cap M:$
for some set $A\in P_I\cap M$, $A\leq b$ and $r\in A\}$. A priori, this set does not have to be a filter on $\ba\cap M,$
but in the situations discussed below it will be. Observe that if $N$ is a model such that $N\cap P_I=M\cap P_I$ and $N\cap\ba$ is dense in $M\cap\ba$
and the real $r$ is $N$-generic, then it is $M,\ba$-weakly generic and the set $\dot H\cap M/r$ is a filter on $\ba\cap M,$
even though not necessarily an $M$-generic filter.

The key claim:

\proclaim {1.5. Claim}
There is an $\in$-sequence $\langle M_k:k\in\omega\rangle$ of countable elementary submodels of a sufficiently large structure such that
for every infinite set $x\subset\omega$ the following set $A_x$ is $I$-positive: $A_x=\{r\in\R:$ for all $k\in\omega$ the real is 
$M_{k+1},\ba$-weakly generic,
and $k\in x\leftrightarrow$
the set $\dot H\cap M_{k+1}/r$ is an $M_{k+1}$-generic filter on $\ba\cap M_{n+1}\}$.
\endproclaim

Suppose that the claim has been proved and $M_k:k\in\omega$ are the ascertained models and $M$ their union.
Then for distinct infinite sets $x,y\subset\omega,$
the sets $A_x$ and $A_y$ are disjoint and even more than that, if $r\in A_x$ and $s\in A_y$ are reals then the filters $\dot H\cap M/r$ and
$\dot H\cap M/s$ are distinct subsets of the poset $B\cap M.$ I am going to find a $I$-positive Borel set $B\subset\R$ such that
$B\subset\bigcup_xA_x$ and for every infinite set $x\subset\omega$ the intersection $B\cap A_x$ contains at most one element. This will complete
the proof since by an absoluteness argument between $V$ and $V[G]$, $B\Vdash\dot r_{gen}$ is 
the unique real $r\in\dot B$ such that $\dot H\cap\check M=
(\dot H\cap M)/r.$ By another absoluteness argument between $V[H]$ and $V[G]$, this unique real must belong to the model $V[H]$.
In other words $B$ forces that $\dot r_{gen}$ can be reconstructed from $\dot H\cap\check M$ and so $V[G]=V[H].$

To find the set $B\subset\R,$ note that the relation
$r\in A_x$ is Borel.
Let $U\subset[\omega]^{\aleph_0}\times\R$ be a $\Sigma^1_n$ universal set and using the projective uniformization find a projective function $f$ such that
$\dom(f)=\{x\subset\omega:$ the vertical section $U_x$ of $U$ belongs to the ideal $I\}$ and for each $x\in\dom(f)$ the value
$f(x)$ is an element of the $I$-positive set $A_x\setminus U_x.$ Look at the projective set rng$(f)\subset\R.$ This set must be $I$-positive
since every generator of the ideal $I$ is of the form $U_x$ for some infinite set $x\subset\omega$ and then $f(x)\in\rng(f)$
is a real that does not belong to that generator. Now just let $B\subset\R$ be any Borel $I$-positive subset of $\rng(f)$ using the assumption (2).

All that remains to be done is the verification of the claim. First construct an $\in$-tower $\langle N_\alpha:\alpha\in\omega_1\rangle$
of countable elementary submodels of a large enough structure so that

\roster
\item at successor ordinals $\alpha$ $\langle N_\beta:\beta\in\alpha\rangle\in N_\alpha$ and 
at limit ordinals $\alpha$ $\omega_1\cap N_\alpha=\bigcup_{\beta\in\alpha}(\omega_1\cap N_\beta)$.
\item if $\alpha=\omega_1\cap N_\alpha$ then whenever possible subject to (1) the model $N_\alpha$ is such that there is a sequence 
$\langle N^\prime_\beta:\beta\in\omega_1\rangle$ in the model $N_\alpha$ with $N_\beta=N^\prime_\beta$ for all $\beta\in\alpha.$
The tower will necessarily be discontinuous at such ordinals.
\item at limit ordinals $\alpha$ where (2) does not happen $N_\alpha=\bigcup_{\beta\in\alpha}N_\beta.$
\endroster

Let me denote the set of all points $\alpha\in\omega_1$ at which (2) happens by $D.$ The set $D$ is uncountable. For if it were not, consider
a countable elementary submodel $M$ of the structure containing the sequence  $\langle N_\alpha:\alpha\in\omega_1\rangle$.
Letting $\alpha=\omega_1\cap M$ it should be that $\alpha$ is greater than all points in $D$ but at the same time, $\alpha\in D$ as witnessed 
by the model $M.$ Contradiction.

Now let $\langle\alpha_k:k\in\omega\rangle$ be the first $\omega$ many ordinals in the set $D,$ and let $M_k=N_{\alpha_k}$ for every number $k\in\omega.$
This is the required sequence of models, but why should the sets $A_x$ be $I$-positive? For any infinite set $x\subset\omega$ consider the
continuous $\in$-tower $T_x$ indexed by the ordinals in the set $\{\alpha_{k+1}:k\in x\}\cup\bigcup\{(\alpha_k,\alpha_{k+1}]:k\notin x\}$.
If $\beta$ is an ordinal in this set then the $\beta$-th model on this tower is just $N_\beta$ unless $\beta=\alpha_{k+1}$ for $k\notin x$, where the
$\beta$-th model is $\bigcup_{\gamma\in\beta} N_\gamma$ as dictated by the continuity requirement. Let $B_x=\{r\in\R:r$ is
generic for every model on the tower $T_x\}$. This set is $I$-positive by the argument from Lemma 1.2. The proof of the claim will be complete once 
I show that $B_x\subset A_x.$

Let $r\in B_x$ be a real. I must verify that $r\in A_x.$ Well, if $k\in x$ then the model $M_{k+1}$ is on the tower $T_x$, the real
$r$ is $M_{k+1}$-generic, therefore $M_{k+1},\ba$-weakly generic and the set $\dot H\cap M_{k+1}/r$
is an $M_{k+1}$-generic filter on $\ba\cap M_{k+1}$ as required in the definition of the set $A_x.$ But what if $k\notin x$?
Look at the model $M_{k+1}$ and choose in it a sequence $\langle N^\prime_\beta:\beta\in\omega_1\rangle$ such that for all $\beta\in
\alpha_{n+1}$ $N_\beta=N^\prime_\beta$ holds. Now since the algebra $\ba$ is nowhere c.c.c. and of density $\aleph_1=2^{\aleph_0}$, it must be that
$\ba\Vdash$ for cofinally many ordinals $\beta\in\omega_1$ the filter $\dot H\cap\check N^\prime_\beta$ is not $\check N^\prime_\beta$-generic.
But for all ordinals
$\beta$ between $\alpha_k$ and $\alpha_{k+1}=\omega_1\cap M_{k+1}$ the models $N'_\beta=N_\beta$ are on the tower $T_x$ and so both the real $r$
and the filter $\dot H\cap N_\beta/r$ are $N_\beta$-generic. This means that (a) the real $r$ is $M_{k+1},\ba$-weakly generic; 
and (b) the filter $\dot H\cap M_{k+1}/r$ cannot be $M_{k+1}$-generic
by the elementarity of the model $M_{k+1}$. Thus $r\in A_x$ as required. \qed
\enddemo

The assumptions of the Lemma feel somewhat ad hoc. I do not have any example of an ideal $I$ satisfying the assumptions that would not be generated
by analytic sets. I do not have an example of a definable ideal $I$ such that the properness of the poset $P_I$ would not be absolute
throughout forcing extensions. I also do not have an example of a definable ideal $I$ such that the forcing $P_I$ is proper but not $<\omega_1$-proper.
The conjecture though is that even in the presence of large cardinals there are such ideals.

The Lemma can be applied to posets like Laver forcing, if there is a suitable determinacy argument that verifies (1) and (2) of the Lemma for the
poset. In the case of Laver forcing this has been done in \cite {Z1, Section 3.2}. There is a fine line dividing the definable forcings into two groups:
The $P_I$'s for simply generated ideals $I$, and $P_I$'s for ideals $I$ for which no generating family consisting of simple sets can be found.

\proclaim {1.6. Example} 
Assume that suitable large cardinals exist. Let $I$ be the ideal of sets of subsets of $\omega$ which are nowhere dense in the algebra
$\P(\omega)$ modulo finite. Then for every $n\in\omega$ there is $m\in\omega$ and a boldface $\Sigma^1_m$
set in $I$ which is not a subset of a $\Sigma^1_n$ set in $I.$ 
\endproclaim

\demo {Proof}
Consider the Mathias forcing. By \cite {Z1, Section 3.4}, this forcing is equivalent to $P_I$, 
and every projective $I$-positive set has a Borel $I$-positive subset.
Also the Mathias forcing is $<\omega_1$-proper. If the statement in 1.6 failed then Lemma 1.4 could be applied to say that all the intermediate
extensions of the Mathias real extension are c.c.c. However, Mathias forcing can be decomposed into an iteration of a $\sigma$ closed and c.c.c.
forcing, and the first step in that iteration is certainly not c.c.c. A contradiction. \qed
\enddemo

To argue for Theorem 0.1, fix a $\sigma$-ideal $I$ $\sigma$-generated by a projective collection
of closed sets. Lemmas 1.2 and 1.3 show that the assumptions of Lemma 1.4 are satisfied and so if $r_{gen}$ is a $V$-generic real for the poset $P_I$
and $s\in V[r_{gen}]$ is an arbitrary real, then $V[s]$ is a c.c.c. extension of $V$ or $V[s]=V[r_{gen}].$ Let us investigate the case
of $V[s]$ being a c.c.c. extension of $V.$ Such a real $s$ is obtained through a ground model $I$-positive Borel set $B$ and a Borel function
$f:B\to\R$ such that $B\Vdash\dot f(\dot r_{gen})=\dot s.$ Move back into the ground model and let $J=\{A\subset\R$ Borel: $B\Vdash\dot s\notin\dot A\}.$
Clearly, $J$ is a $\sigma$-ideal of Borel sets and the poset $P_J$ is c.c.c.: an uncountable antichain in it would give an uncountable
antichain in the algebra generated by the name $\dot s.$ Since $P_J$ is c.c.c. and the real $\dot s$ is forced to fall out of all $J$-small
ground model coded Borel sets, the real $\dot s$ is actually forced to be $P_J$-generic. Now I will show that in a Cohen extension there is
a generic real for the poset $P_J$, which will conclude the argument since all complete subalgebras of the Cohen algebra have countable
density and therefore are Cohen themselves. Let $c$ be a $V$-generic Cohen real. There is
in $V[c]$ a real number $d\in B$ which falls out of all ground moded coded $I$-small sets. To see this, let $T\subset(\omega\times\omega)^{<\omega}$ 
be a tree projecting into the set $B$. By thinning out the set $B$ and pruning the tree $T$ if necessary we may assume that for every node $\tau\in T$ 
the projection of the tree $T\restriction\tau$ is still an $I$-positive set. Then viewing $c$ as a generic path through the tree $T$ an almost trivial
density argument shows that its first coordinate $d$ does fall out of all ground model $I$-small sets. Look at the real $f(d).$ Whenever $A$ is a ground model
coded $J$-small set, the set $f^{-1}A$ is a ground model coded $I$-small set and so $d\notin f^{-1}A$ and $f(d)\notin A.$ Thus the real $f(d)$
falls out of all ground model coded $J$-small sets and must be generic for $P_J$ as required.

The last remark. Turning the history of forcing on its head, the understanding of the forcing $P_I$ means finding a determinacy argument that will produce
a dense subset of $P_I$ consisting of combinatorially manageable sets, for example perfect sets in the case of Sacks forcing. 
Remarkably, in all known cases this also leads to the proof of the following proposition: for every $I$-positive Borel set $B$ there is a Borel
function $f:\R\to B$ such that the preimages of $I$-small sets are $I$-small. This property of the ideal $I$ is critical in the proof that
the covering number for $I$ can be isolated, see \cite {Z1}. Can such a feat be repeated for ideals like the closed measure zero ideal?


\Refs\widestnumber\key{Z2}
\ref
 \key B
 \by T. Bartoszynski and H. Judah
 \book Set Theory: On the Structure of the Real Line
 \yr 1995
 \publ A K Peters
 \publaddr Wellesley, Massachusets
\endref
\ref
 \key J
 \by T. Jech
 \book Set Theory
 \yr 1978
 \publ Academic Press
 \publaddr New York
\endref
\ref
 \key M
 \by D. A. Martin and J. Steel
 \paper A proof of projective determinacy
 \jour J. Amer. Math. Soc.
 \yr 1989
 \vol 85
 \pages 6582--6586
\endref
\ref
 \key N
 \by I. Neeman and J. Zapletal
 \paper Proper forcings and absoluteness in $L(\R)$
 \jour Comment. Math. Univ. Carolinae
 \yr 1998
 \vol 39
 \pages 281--301
\endref
\ref
 \key S
 \by S. Solecki
 \paper Covering analytic sets by families of closed sets
 \jour J. Symbolic Logic
 \vol 59
 \yr 1994
 \pages 1022--1031
\endref
\ref
 \key W
 \by W. H. Woodin
 \paper Supercompact cardinals, sets of reals and weakly
homogeneous trees
 \jour Proc. Natl. Acad. Sci. USA
 \vol 85
 \pages 6587--6591
 \yr 1988
\endref
\ref
 \key Z1
 \by J. Zapletal
 \paper Isolating cardinal invariants
 \jour J. Math. Logic
 \paperinfo accepted
\endref
\ref
 \key Z2
 \bysame
 \paper Countable support iteration revisited
 \jour J. Math. Logic
 \paperinfo submitted
\endref
\endRefs
\enddocument